\documentclass[11pt,english]{smfart}
   
\usepackage{dsfont,eucal,mathptmx,txfonts}

\usepackage[text={6.5in,9.2in},centering]{geometry}

\usepackage[dvipsnames]{xcolor}   
\usepackage{xparse}
\usepackage{xr-hyper}
\definecolor{brightmaroon}{rgb}{0.76, 0.13, 0.28}
\usepackage[linktocpage=true,colorlinks=true,hyperindex,citecolor=teal,linkcolor=blue]{hyperref}   
    
\renewcommand{\labelenumi}{\upshape(\arabic{enumi})}

\newtheorem{thm}{Theorem}[section]

\newtheorem{prop}[thm]{Proposition}
 
\newtheorem{lem}[thm]{Lemma}

\numberwithin{equation}{section}
   
\author{Amartya Goswami} 
\address{[1] Department of Mathematics and Applied Mathematics,
University of Johannesburg, 
P.O. Box 524, 
Auckland Park 
2006,
South Africa \\ [2] National Institute for Theoretical and Computational Sciences (NITheCS), South Africa.}

\email{agoswami@uj.ac.za}  
\title{Proper spaces are spectral}   



\begin{document} 

\begin{abstract}
Since Hochster's work, spectral spaces have attracted increasing interest.
Through this note we give a new self-contained and constructible topology-independent proof
of the fact that the set of proper ideals of a ring endowed with coarse lower topology is a spectral space. 
\end{abstract}   
\maketitle 
 
\section{Introduction}
  
Hochster \cite{H69} has characterised spectral spaces (\textit{i.e.}, quasi-compact, sober,  admitting a basis of quasi-compact open subspaces that is closed under finite intersections) as spectra of prime ideals of rings endowed with Zariski topology. Since then  spectrality has been  studied in different contexts. Here are some examples. Apart from Zariski  topology, several other topologies (like Lawson, Scott, lower, upper, ultra-filter, patch, inverse \textit{etc}.) have been considered on prime ideals (see \cite{DST19} for a detailed account). Among other classes of ideals, the set of all ideals with lower topology is a spectral space has been proved in \cite{P94},  whereas \cite{FFS2016} has shown that the spaces of all ideals, proper ideals, and radical ideals of a ring endowed with a hull-kernel-type topology are spectral. Recall also that the Zariski-Riemann space of valuation domains of a given field is one of the first examples of a ``special'' spectral space (see \cite{DF86, DFF87}). In \cite{FFS16}, the space of the semistar operations (of finite type), endowed with a Zariski-like   topology has been shown to be spectral.
In \cite{R20}, it has been shown that the collection of continuous valuations on a topological monoid with topology determined by any finitely generated ideal is a spectral space.
In \cite{AH12}, the authors have considered the prime spectra of modules with different
topologies and introduced some conditions under which these are spectral spaces. Another context of spectral spaces is Stone duality: the prime filters of a distributive lattice with the Stone topology is a spectral space (see \cite[Chapter 3]{DST19})   and the proper filters of a Boolean algebra (or more generally an
ortholattice) with a Stone-like topology is a spectral space (see \cite{BH20, MY22}).

Our aim is to prove that the set of proper ideals of a ring endowed with a coarse lower topology is spectral. Our proof is self-contained and constructible topology-independent.To prove the spectrality, we have adapted the technique that avoids the  checking of the existence of a basis of quasi-compact open subspaces that is closed under finite intersections. These properties rather follow from proper ideals to be an open subspace of the spectral space of all  ideals (see the discussion before Lemma \ref{cso}). 

Note that this spectrality can also be deduced by combining some propositions contained in \cite{FFS16} and in \cite{FS20}, where it is observed that the inverse topology of the Zariski topology of $\mathrm{Idl}(A)$ is precisely the coarse lower topology, therefore both these topologies have the same constructible topology. A coarse lower topology on an arbitrary class of ideals of a ring has been studied in \cite{DG22}, and these classes of spaces are called ideal spaces. Although some of the results appearing here have been discussed in \cite{DG22}, to make the present work self-contained, we have provided proofs of those results as well. 

\section{Proper spaces and spectrality}

Suppose $A$ is a commutative ring with identity and $\mathrm{Idl}(A)$ denotes the set of all ideals of $A$. For $\mathfrak{I}\in \mathrm{Idl}(A),$ we define $\mathrm{v}(\mathfrak{I})=\{\mathfrak{S}\in \mathrm{Idl}(A) \mid \mathfrak{I}\subseteq \mathfrak{S}\}.$ The collection $\mathcal{S}_{\mathrm{v}}=\{\mathrm{v}(\mathfrak{I})\mid\mathfrak{I}\in \mathrm{Idl}(A)\}$ of subsets of $\mathrm{Idl}(A),$ as a closed subbase, induces a unique (coarse lower) topology on $\mathrm{Idl}(A).$ The set $\mathrm{Idl}(A)$ endowed with the coarse lower topology is called an \emph{ideal space}. Let $\mathrm{Prp}(A)$ be the set of all proper ideals (\textit{i.e.}, $A\notin \mathrm{Prp}(A)$) of $A$. The set $\mathrm{Prp}(A)$ with the subspace topology of the ideal space $\mathrm{Idl}(A)$ is called a \emph{proper space}. 
It is well known (see \cite[Theorem 4.2]{P94}) that 
  
\begin{thm}\label{idsp} 
Every ideal space $\mathrm{Idl}(A)$ is  spectral.
\end{thm}
 
The proof of Theorem \ref{idsp} is based on the fact that $\mathrm{Idl}(A)$ is an algebraic lattice, and that is no longer true for $\mathrm{Prp}(A)$.
The proof of spectrality of $\mathrm{Prp}(A)$ relies on the following key Lemma \ref{cso}, which allows us to use the spectrality of $\mathrm{Idl}(A)$. Furthermore, we can avoid proving the existence of a compact open base and closedness of open compact sets under finite intersections. Thanks to Lemma \ref{cso}, these two properties follow from the spectrality of  $\mathrm{Idl}(A).$ 

\begin{lem}\label{cso}
A quasi-compact, sober, open subspace of a spectral space is spectral. 
\end{lem}

\begin{proof}
Let $S$ be a quasi-compact, sober, open subspace of a spectral space $X$. Since $S$ is quasi-compact and sober, all we need to prove is that the set $\mathcal{B}_{\scriptscriptstyle S}$ of compact open subsets of $S$ forms a basis of a topology that is closed under finite intersections. The following two facts are obvious:
\begin{itemize}
\item[$\bullet$] A subset of $S$ is open in $S$ if and only if it is open in $X$.
\item[$\bullet$] Therefore, a subset of $S$ belongs to $\mathcal{B}_{\scriptscriptstyle S}$ if and only if it belongs to $\mathcal{B}_{\scriptscriptstyle X}.$
\end{itemize} 
Using these facts, we argue as follows:
 
Let $U$ be an open subset of $S$. Since $U$ is also open in $X$, we have $U=\cup\, \mathcal{U},$ for some subset $\mathcal{U}$ of $\mathcal{B}_{\scriptscriptstyle X}.$ But each element of $\mathcal{U}$ being a subset of $U$ is a subset of $S$, and it belongs to $\mathcal{B}_{\scriptscriptstyle S}.$ Therefore every open subset of $S$ can be presented as a union of compact open subsets of $S$. Now it remians to prove that $\mathcal{B}_{\scriptscriptstyle S}$ is closed under finite intersections, but this immediately follows from the fact that $\mathcal{B}_{\scriptscriptstyle X}$ is closed under finite intersections. 
\end{proof} 

To prove the claim that a proper space $\mathrm{Prp}(A)$ is spectral, now it is sufficient to show that:

\begin{itemize}
\item[$\mathrm{(I)}$] \label{oco} $\mathrm{Prp}(A)$ is quasi-compact;
\item[$\mathrm{(II)}$] \label{tso} $\mathrm{Prp}(A)$ is sober;
\item[$\mathrm{(III)}$] \label{top}  $\mathrm{Prp}(A)$ is an open subspace of the ideal space $\mathrm{Idl}(A)$.
\end{itemize}

(I).  
Notice that $\mathrm{Prp}(A)$ satisfies the following property for all subbasic closed sets $\mathrm{v}(\mathfrak{I})$:
\[\mathfrak{I}\in \mathrm{Idl}(A),\,\mathrm{v}(\mathfrak{I})=\emptyset\,\,  \text{implies}\,\, \mathfrak{I}=A,\label{star} \tag{$\star$}\]
which we need in order to apply the finite intersection property. Since the aurgument goes with subbasic closed sets, we also require the Alexander subbase theorem here.

\begin{prop}\label{comp} 
Every proper space $\mathrm{Prp}(A)$ is quasi-compact. 
\end{prop} 

\begin{proof}   
Let  $\{C_{\scriptscriptstyle  \alpha}\}_{\scriptscriptstyle \alpha \in \Lambda}$ be a family of subbasic closed sets of an ideal space $\mathrm{Prp}(A)$   such that $\bigcap_{\scriptscriptstyle \alpha\in \Lambda}C_{ \scriptscriptstyle \alpha}=\emptyset.$ Let $\{\mathfrak{I}_{ \scriptscriptstyle \alpha}\}_{\scriptscriptstyle \alpha \in \Lambda}$ be a family of ideals of $\mathrm{Idl}(A)$ such  that $\forall \alpha \in \Lambda,$  $C_{\scriptscriptstyle  \alpha}=\mathrm{v}(\mathfrak{I}_{\scriptscriptstyle  \alpha}).$  Since $\bigcap_{\alpha \in \Lambda}\mathrm{v}(\mathfrak{I}_{\scriptscriptstyle  \alpha})=\mathrm{v}\left(\sum_{\scriptscriptstyle \alpha \in \Lambda}\mathfrak{I}_{\scriptscriptstyle  \alpha}\right),$ we get  $\mathrm{v}\left(\sum_{\scriptscriptstyle \alpha \in \Lambda}\mathfrak{I}_{\scriptscriptstyle  \alpha}\right)=\emptyset,$ and that (\ref{star}) gives $ \sum_{\scriptscriptstyle \alpha \in \Lambda}\mathfrak{I}_{\scriptscriptstyle  \alpha}=A.$ Then, in particular, we obtain $1=\sum_{\scriptscriptstyle \alpha_i\in \Lambda}x_{\scriptscriptstyle  \alpha_i},$ where $x_{ \scriptscriptstyle \alpha_i}\in \mathfrak{I}_{\scriptscriptstyle \alpha_i}$ and $x_{\scriptscriptstyle  \alpha_i}\neq 0$ for $i=1, \ldots, n$. This implies    $A=\sum_{\scriptscriptstyle   i \, =1}^{\scriptscriptstyle  n}\mathfrak{I}_{\scriptscriptstyle \alpha_i}.$ Therefore,   $\bigcap_{\scriptscriptstyle  i\,=1}^{\scriptscriptstyle  n}C_{\scriptscriptstyle  \alpha_i}=\emptyset,$ and hence by the Alexander subbase theorem, $\mathrm{Prp}(A)$ is quasi-compact.  
\end{proof}   
  
(II). 
To show $\mathrm{Prp}(A)$ is sober, we show that $\mathrm{Prp}(A)$ is $T_{\scriptscriptstyle 0}$ and every non-empty irreducible closed subset of $\mathrm{Prp}(A)$ has a generic point.   Since  for a topological space to be $T_{\scriptscriptstyle 0}$, it is necessary and sufficient that its specialisation order is a partial order (see \cite{P94}), a proper space $\mathrm{Prp}(A)$ has the same property. In other words,
 
\begin{prop}
Every proper space $\mathrm{Prp}(A)$ is   $T_{\scriptscriptstyle 0}.$ 
\label{ct0t1}  
\end{prop} 

For the proof of existence of generic points of non-empty irreducible closed sets, we need the following result.

\begin{lem}\label{irrc}
Every non-empty subbasic closed subset of $\mathrm{Prp}(A)$ is irreducible. 
\end{lem} 
 
\begin{proof}  
We show that $\mathrm{v}(\mathfrak{I})=\mathrm{cl}(\mathfrak{I})$ for every non-empty subbasic closed set $\mathrm{v}(\mathfrak{I})$ of $\mathrm{Prp}(A),$ and that will prove the lemma. Since $\mathrm{cl}(\mathfrak{I})$ is the smallest closed set containing $\mathfrak{I}$, and since $\mathrm{v}(\mathfrak{I})$ is a closed set containing $\mathfrak{I}$, we immediately have  $\mathrm{cl}(\mathfrak{I})\subseteq \mathrm{v}(\mathfrak{I})$. 
For the reverse inclusion, we first consider the case that $\mathrm{cl}(\mathfrak{I})= \mathrm{Prp}(A)$. Then  the above containment forces $\mathrm{v}(\mathfrak{I})=\mathrm{cl}(\mathfrak{I})$. Now suppose that $\mathrm{cl}(\mathfrak{I})\neq \mathrm{Prp}(A)$. Being a closed set, $\mathrm{cl}(\mathfrak{I})$ is expressible as   
$
\mathrm{cl}(\mathfrak{I})={\bigcap_{\scriptscriptstyle\alpha\in\Omega}}\left({\bigcup_{\scriptscriptstyle i\,=1}^{\scriptscriptstyle n_\alpha}}\mathrm{v}(\mathfrak{I}_{\scriptscriptstyle\alpha i})\right),
$
for some index set $\Omega$, such that,  for each $\alpha\in\Omega$, there is a positive integer $n_{\alpha}$ and ideals $\mathfrak{I}_{\scriptscriptstyle \alpha 1},\dots, \mathfrak{I}_{\scriptscriptstyle \alpha n_\alpha}$ of $\mathrm{Idl}(A).$ 
Since 
$\mathrm{cl}(\mathfrak{I})\neq \mathrm{Prp}(A)$, without loss of generality, we assume that for each $\alpha$, ${\bigcup_{ \scriptscriptstyle i\,=1}^{\scriptscriptstyle  n_\alpha}}\mathrm{v}(\mathfrak{I}_{\scriptscriptstyle \alpha i})$ is non-empty. Therefore $\mathfrak{I}\in   {\bigcup_{\scriptscriptstyle  i\,=1}^{\scriptscriptstyle  n_\alpha}}\mathrm{v}(\mathfrak{I}_{\scriptscriptstyle \alpha i})$ for each $\alpha$, from which we have $\mathrm{v}(\mathfrak{I})\subseteq {\bigcup_{\scriptscriptstyle  i=1}^{\scriptscriptstyle  n_\alpha}}\mathrm{v}(\mathfrak{I}_{\scriptscriptstyle \alpha i})$ for each $\alpha$. This implies $\mathrm{v}(\mathfrak{I})\subseteq \mathrm{cl}(\mathfrak{I})$. 
\end{proof}    

If $C$ is an irreducible closed subset of  a topological space $X$ and $\mathcal{S}$ is a closed subbase of $X$, then it is known (see \cite[Section 7.2, p.\,308]{H73}) that $C$ is the intersection of members of $\mathcal{S}.$ For a proper space we get more in the sense that the converse of Lemma \ref{irrc} is also true. 

\begin{lem}\label{ircs} 
If $C$ is a non-empty irreducible closed subset of a proper space $\mathrm{Prp}(A)$, then $C=\mathrm{v}(\mathfrak{I})$ for some $\mathfrak{I}\in \mathrm{Prp}(A).$ 
\end{lem}

We are now ready to prove the sobriety of $\mathrm{Prp}(A)$.
 
\begin{prop}\label{sob}  
Every proper space $\mathrm{Prp}(A)$ is sober.
\end{prop}  

\begin{proof}  
Since by Lemma \ref{ircs}, every non-empty irreducible closed subset of $\mathrm{Prp}(A)$ is of the form $\mathrm{v}(\mathfrak{I})$, where $\mathfrak{I}\in \mathrm{Prp}(A).$ Suppose $\mathrm{v}(\mathfrak{I})$ be a non-empty irreducible closed subset of $\mathrm{Prp}(A).$ From the fact $\mathfrak{I}\in\mathrm{v}(\mathfrak{I})$, we immediately have $\mathrm{cl}(\mathfrak{I})\subseteq \mathrm{v}(\mathfrak{I}).$ Therefore, to show $\mathrm{v}(\mathfrak{I})$ has a generic point, it is sufficient to prove that $\mathrm{cl}(\mathfrak{I})\supseteq \mathrm{v}(\mathfrak{I}).$ Since $\mathcal{S}_{\mathrm{v}}$ is a closed subbase of $\mathrm{Prp}(A)$, the required containment now follows from Lemma  \ref{irrc}. Moreover, by Proposition \ref{ct0t1}, an ideal space $\mathrm{Prp}(A)$ is $T_{\scriptscriptstyle  0}$. Therefore, we obtain the uniqueness of a  generic point. 
\end{proof}  

(III). What remains is  to show that  
 
\begin{prop}
\label{open}
A proper space $\mathrm{Prp}(A)$ is open in $\mathrm{Idl}(A)$.
\end{prop}

\begin{proof} 
Since $A\in\mathrm{Idl}(A),$ by Lemma \ref{irrc}, $A=\mathrm{v}(A)=\mathrm{cl}(A),$ and therefore $\mathrm{Idl}(A) {\setminus}\mathrm{Prp}(A)$ is closed, and that implies $\mathrm{Prp}(A)$ is open. 
\end{proof}

Since $\mathrm{Prp}(A)$ is quasi-compact (Proposition \ref{comp}), sober (Proposition \ref{sob}), and an open subspace of $\mathrm{Idl}(A)$ (Proposition \ref{open}), by Lemma \ref{cso}, we have 

\begin{thm}
Every proper space is spectral.
\end{thm} 

\section*{Acknowledgement} 
The author is indebted to \textsc{George Janelidze} for some fruitful discussions.


\begin{thebibliography}{1} 
\bibitem{AH12} \textsc{A. Abbasi} and \textsc{D. Hassanzadeh-Lelekaami}, Modules and spectral spaces, 
\textit{Comm. Algebra}, \textbf{40}(11) (2012), 4111--4129.

\bibitem{BH20}
\textsc{N. Bezhanishvili} and \textsc{W. H. Holliday}, Choice-free Stone duality, \textit{J. Symb. Log.}, \textbf{85} (2020), 109--148, 2020.

\bibitem{MY22}
\textsc{J. McDonald} and \textsc{K. Yamamoto},
Choice-free duality for orthocomplemented lattices by means of spectral spaces, Algebra Universalis, \textbf{83}(37) (2022), paper no. 37, 32 pp.

\bibitem{DF86}
\textsc{D. E. Dobbs} and \textsc{M. Fontana} Kronecker function rings and abstract Riemann surfaces, \textit{J. Algebra}, \textbf{99} (1986), 263--274.

\bibitem{DFF87}
\textsc{D. E. Dobbs}, \textsc{R. Fedder}, and \textsc{M. Fontana}, Abstract Riemann surfaces of integral domains and spectral spaces, \textit{Annali Mat. Pura Appl.}, \textbf{148} (1987), 101--115.

\bibitem{DST19}  \textsc{M. Dickmann}, \textsc{N. Schwartz}, and \textsc{M. Tressel},  \emph{Spectral spaces}, Cambridge Univ. Press, 2019. 

\bibitem{FFS2016} \textsc{C. A. Finocchiaro}, \textsc{M. Fontana}, and \textsc{D. Spirito},  A topological version of Hilbert's Nullstellensatz, \textit{J. Algebra}, \textbf{461} (2016), 25--41. 

\bibitem{FFS16} \textsc{C. A. Finocchiaro}, \textsc{M. Fontana}, and \textsc{D. Spirito}, New distinguished classes of spectral spaces: a survey, in: \textit{Multiplicative ideal theory and factorization theory}, 117--143,
Springer, 2016.

\bibitem{FS20}
\textsc{C. Finocchiaro} and \textsc{D. Spirito}, Suprema in spectral spaces and the constructible closure, \textit{New York J. Math.}, \textbf{26} (2020), 1064--1092.

\bibitem{DG22} \textsc{T. Dube} and \textsc{A. Goswami}, Ideal spaces: an extension of structure spaces of a ring, \textit{J. Algebra Appl.} (online ready).  
        
\bibitem{H73} \textsc{D. Harris}, Universal quasi-compact $T_{\scriptscriptstyle 1}$ spaces,       
\textit{General Topology and Appl.}, \textbf{3} (1973), 291--318.   
 
\bibitem{H69} 
\textsc{M. Hochster}, Prime ideal structure in commutative rings, \textit{Trans. Am. Math. Soc.}, \textbf{142} (1969), 43--60.
  
\bibitem{P94} \textsc{H. A. Priestley}, Intrinsic spectral topologies, in: \textit{Papers on general topology and applications}
(Flushing, NY, 1992), \textbf{728},  78--95, New York Acad. Sci., New York, 1994.

\bibitem{R20} \textsc{S. Ray}, Closure operations, continuous valuations on monoids and spectral spaces, \textit{J. Algebra Appl.}, \textbf{19}(1) (2020),  2050006, 27 pp.
\end{thebibliography}
\end{document}